\documentclass[english, 12pt]{article}
\usepackage[english]{babel,varioref}
\usepackage[ansinew]{inputenc}
\usepackage[fleqn,reqno]{amsmath}
\usepackage{amsthm}
\usepackage{paralist}
\usepackage{natbib}
\usepackage{graphics}
\usepackage{epic,eepic}

\bibpunct[: ]{(}{)}{,}{a}{}{,}

 \setdefaultenum{(i)}{a)}{(1)}{(A)}

\theoremstyle{plain}

\newtheorem{thm}{Theorem}
\newtheorem{cor}[thm]{Corollary}
\newtheorem{lem}[thm]{Lemma}

\newtheorem{defn}[thm]{Definition}

\def\ld{\mathop{\rm ld}}
\begin{document}

\title{Beyond the law of large numbers: Introducing progressive sampling, weaving, the geometric triangle, and corresponding distributions}

\author{Uwe Saint-Mont}

\maketitle

\vspace{4ex}

{\bf AMS classifiction}: 60E05, 60F99, 62D99

\section*{Abstract}
In probability theory and statistics, the IID model represents a single population, and a large, potentially infinite sample from this population. Main theorems, in particular the central limit theorem and laws of large number (LLN) assure convergence, making asymptotic statistics possible.

To {\it avoid} convergence, it is thus straightforward to consider two populations and a sample that ceaselessly fluctuates between them. It is the aim of this contribution to study the effects that thus occur. To this end, we introduce ``progressive sampling,'' leading to a straightforward model that is analytically tractable.
With a minimum of technical overhead, a number of interesting results thus ensue:

In particular, one encounters a multiplicate structure (similar to Pascal's triangle) that is associated with a new class of distributions (related to the binomial).
Although the argument is completely probabilistic, it entails a well-known fractal structure. It also turns out that the new (global) operation of ``weaving'' is equivalent to a certain (local) cascade process.

\newpage

\section*{Introduction}
Traditional statistics rests on several main theorems, in particular the central limit theorem and laws of large number (LLN). Given an IID sequence $X_1,X_2,\ldots$ of random variables, the basis of Frequentist statistics is the convergence of ${\bar X}_n=S_n /n$,  where $S_n=\sum_{i=1}^n X_i$. However, in calculus, convergence of a sequence $x_1,x_2,\ldots$ is a strong assumption, and, typically, not even the (much weaker) Cesaro-limit $\lim_{n \rightarrow \infty}{\bar x}_n= \lim_{n \rightarrow \infty}(\sum x_i/n)$ exists. In dynamic system theory, also, convergence towards a single point is a rare exception.

In probability theory, the IID model represents a single population, and a large, potentially infinite sample from this population. To {\it avoid} convergence, it is thus straightforward to consider two populations (distributions), $H_0$ and $H_1$ say, and a sample that ceaselessly fluctuates between them. In other words, if one switches between the populations, ${\bar X}_n$ will not converge. (In the jargon of dynamic system theory, the (unique) limit is replaced by a simple attractor.)

However, a constant switching rate won't do: If $j$ observations from $H_0$ are followed by $j$ observations from $H_1$, and so forth, the arithmetic mean of this sequence will converge. However, if $2^j$ observations from $H_0$ are followed by $2^{j+1}$ observations from $H_1$ etc. one has the desired effect. (On a logarithmic scale, taking $\ld =\log_{2}$, the distance $(j+1)-j=1$ is a constant. Thus, there, one switches at a constant rate.)

Altogether one obtains a stochastic process that is inhomogeneous in a particular way. Its paths depend on the concrete distributions $H_0, H_1$ and the way switching is done. The aim of this article is to explore straightforward consequences of this setting.

\section{The weaver's distribution}

In order to keep things finite, suppose for the rest of this contribution that first moments exist, such that without real loss of generality $\mu(H_0)=0$ and $\mu(H_1)=1$ are the expected values of the two distributions involved.

A particularly simple way to alternate between $H_0$ and $H_1$ is to take the next batch of $2^j$ observations $(j=0,1,\ldots)$ from population $H_0$ with probability $1-p$ and from population $H_1$ with probability $p$. Thus one creates a hierarchical random system (a particular random probability measure), composed of a choice mechanism, selecting the population in charge, and the realization mechanism, providing observations from the population selected.

With probability $p$, the first observation comes from $H_1$, and with probability $1-p$, the first observation comes from $H_0$. Thus, conditional on this choice, the expected value observed is either $\mu(H_1)=1$ or $\mu(H_0)=0$, and the unconditional mean is $\mu=p \mu(H_1) + (1-p) \mu(H_0) =p$.

With probability $p$, the second {\it and} the third observations come from $H_1$, and with probability $1-p$, these observations come from $H_0$. Thus, after two choices, the overall situation is as follows:

\begin{tabular}{|l|l|l|l|}
  \hline
   Number of observations & Number of observations & Probability & Conditional \\
   from $H_0$ & from $H_1$ & & Mean \\\hline
    1+2 & 0 & $(1-p)^2$ & 0 \\
    1 & 2 & (1-p) p & $2/3$ \\
    2 & 1 & p (1-p) & $1/3$ \\
    0 & 3 & $p^2$ & 1 \\
  \hline
\end{tabular}

 The unconditional mean does not change, since
$$
\mu = p^2 + \frac{1}{3}p(1-p) + \frac{2}{3} (1-p)p = p^2 +p (1-p) =p.
$$

Notice that similar to the binomial distribution, every path splits in two. However, unlike the binomial distribution, the paths do not combine. Rather, like threads, they interweave.

\vspace{1ex}
{\bf Illustration:}\label{illu} Binomial structure, lokal splitting, and global weaving
$$
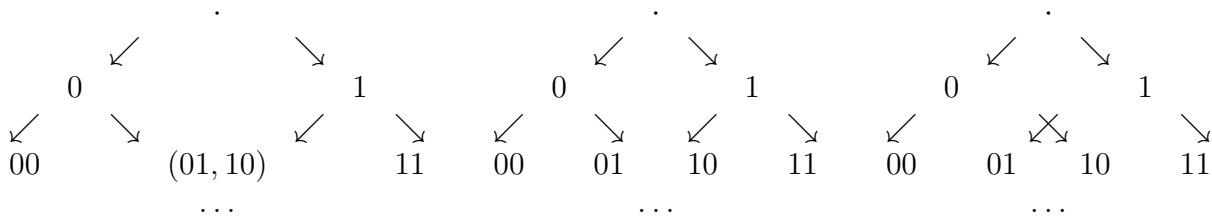
\begin{array}{ccccccc}
&&& \cdot &&& \\
 &  & \swarrow && \searrow &  & \\
 & 0 &&&& 1 & \\
  \swarrow &  & \searrow &  & \swarrow & & \searrow \\
  00 &  &  & (01, 10) &  &  & 11 \\
  \multicolumn{7}{c}{\ldots}
\end{array}
\hspace{3ex}
\begin{array}{ccccccc}
&&& \cdot &&& \\
 &  & \swarrow && \searrow &  & \\
 & 0 &&&& 1 & \\
  \swarrow &  & \searrow &  & \swarrow & & \searrow \\
  00 &  & 01 &  & 10 &  & 11 \\
  \multicolumn{7}{c}{\ldots}
\end{array}
\hspace{3ex}
\begin{array}{ccccccc}
&&& \cdot &&& \\
 &  & \swarrow && \searrow &  & \\
 & 0 &&&& 1 & \\
  \swarrow & \multicolumn{5}{c}{\searrow \hspace{-3ex} \swarrow}  & \searrow \\
  00 &  & 01 &  & 10 &  & 11 \\
    \multicolumn{7}{c}{\ldots}
\end{array}
$$

After $n$ steps (selections, choices), one thus obtains an interesting distribution, called the ``weaver's distribution'' $W(n,p)$, in the sequel.

\begin{defn}\label{progsamp}
  (Progressive sampling)

  Given two distributions $H_0$ and $H_1$, define progressive sampling as follows: A sample of size $2^n -1$, i.e., $X_1; X_2, X_3; X_4, X_5, X_6 , X_7;\ldots; X_{2^{n-1}},\ldots,X_{2^n-1}$, consists of $n$ sub-samples, where the $i$th sub-sample $X_{2^{i-1}},\ldots, X_{2^{i}-1}$ has size $2^{i-1}$ for $i=1,\ldots,n$.

The selection mechanism $B$ chooses $H_1$ with probability $p$, and $H_0$ with probability $1-p$ (independent of anything else). Thus, with these probabilities, the $i$th sub-sample comes from $H_1$ or $H_0$, respectively. Finally, denote by $B_n$ the collection of $n$ such independent choices.
\end{defn}

\begin{thm}\label{weaver}
              (The weaver's distribution)

 Given the situation described in definition \ref{progsamp}, suppose the first moments are $\mu(H_0)=0$ and $\mu(H_1)=1$, respectively.  

For $n=1,2,\ldots$ let $S_n=\sum_{i=1}^{2^n-1} X_i$, ${\bar X}_n=S_n / (2^n-1)$, and $Y_n = E({\bar X}_n | B_n)$. Some elementary properties of these processes are:
\begin{enumerate}
  \item
  $Y_n$ assumes the values $y_k=y_{k,n}=k/(2^n-1)$ for $k=0,1,\ldots,2^n-1$, and
the difference between the realizations of $Y_n$ is a constant; more precisely,

 $y_{k+1}-y_k = \frac{k+1}{2^n-1}-\frac{k}{2^n-1} = 1/(2^n-1)$ for $k=0,\ldots,2^n-2$

  \item Suppose $B_n={\bf b}_n$, then ${\bf b}_n=(b_{n-1},\ldots,b_1,b_0)$ is a binary vector of length $n$, i.e., $b_{i-1}=0$ if in the $i$th selection, $B$ chooses $H_0$, and $b_{i-1}=1$ otherwise. Note that $b_{i-1}$ can also be interpreted as the $i$th digit in the binary representation of a natural number $k \in \{0,\ldots,2^n-1\}$, i.e., $k=\sum_{i=0}^{n-1} b_i 2^i$. Then the probability $p_k$ at the point $y_k$ is given by
  $$
  p_k = p^{\#1} (1-p)^{\#0} =  p^{\sum_{i=0}^{n-1} b_i}(1-p)^{n-\sum_{i=0}^{n-1} b_i} \ge 0 ,
  $$
where $\#1$ and $\#0$ denote the number of ones and zeros in ${\bf b}_n$, respectively. In particular, every $p_k$ can be written in the form $p_k=p^j (1-p)^{n-j}$ with some $j\in \{0,\ldots,n\}$.

  \item More generally and explicitly, the distributions of $B_n$, $E(S_n|B_n)$, and ${Y}_n$ are
$$
\begin{array}{|r|r|r|c|c|c|}\hline
 (k)_{10} & (k)_2 & {\bf b}_n & E(S_n|{\bf b}_n) & y_{k,n} & p_k \\\hline
 0 & 0 & (0,\ldots,0) & 0 & 0 & (1-p)^n \\
 1 & 1 & (0,\ldots,0,1) & 1 & 1/(2^n-1) & p(1-p)^{n-1} \\
 2 & 10 & (0,\ldots,0,1,0) & 2 & 2/(2^n-1) & p(1-p)^{n-1} \\
 3 & 11 & (0,\ldots,0,1,1) & 3 & 3/(2^n-1) & p^2(1-p)^{n-2} \\
 4 & 100 & (0,\ldots,0,1,0,0) & 4 & 4/(2^n-1) & p(1-p)^{n-1} \\
 {\ldots}&{\ldots}&{\ldots}&{\ldots}&{\ldots}&{\ldots} \\
 2^n-5 & 1\ldots 1011 & (1,\ldots,1,0,1,1) & 2^n-5 & (2^n-5)/(2^n-1) & p^{n-1}(1-p) \\
 2^n-4 & 1\ldots 100 & (1,\ldots,1,0,0) & 2^n-4 & (2^n-4)/(2^n-1) & p^{n-2}(1-p)^2 \\
 2^n-3 & 1\ldots 101 & (1,\ldots,1,0,1) & 2^n-3 & (2^n-3)/(2^n-1) &  p^{n-1}(1-p) \\
 2^n-2 & 1\ldots 10 & (1,\ldots,1,0) & 2^n-2 & (2^n-2)/(2^n-1) & p^{n-1}(1-p) \\
 2^n-1 & \underbrace{1\ldots 1} & (1,\ldots,1) & 2^n-1 & 1 & p^{n} \\
 & n \; \mbox{times} &&&& \\\hline
\end{array}
$$

\end{enumerate}
  \end{thm}

Proof: (i) is obvious since $E(S_n|B_n)$ assumes the values $0,1,\ldots,2^n-1$, and (ii) follows from (iii). (iii) holds by construction, or since by the Binomial theorem
$
\sum_{k=0}^{2^n-1} p_k = \sum_{j=0}^n \binom{n}{j} p^j (1-p)^{n-j} =1.
$
$\diamondsuit$

\vspace{2ex}
We say that $Y_n$ has a {\it weaver's distribution}, $Y_n \sim W(n,p)$, with parameters $n$ and $p$. Since powers of two play a major role, ``binary distribution'' would also be a suitable choice - much in line with ``Bernoulli'' and ``Binomial'', distributions, that are closely related.

\vspace{2ex}
\begin{thm}\label{pascal}
(The geometric triangle)

Given the assumptions and the notation of the last theorem, let ${\bf b}_n=s_{ij}$ be a vector with exactly $i$ ones and $j$ zeros, such that $i+j=n$. Moreover, set $f= p / (1-p)$.
\begin{enumerate}
  \item  The probabilities $p(\cdot)$ of the concatenated vectors $(s_{ij},1), (1,s_{ij}), (s_{ij},0)$, and $(0,s_{ij})$ are: $$
     \frac{p(s_{ij},1) }{ p(s_{ij},0)}=\frac{p(1,s_{ij})}{p(0,s_{ij})} = \frac{ p^{i+1}(1-p)^j }{ p^i(1-p)^{j+1}} = \frac{p}{1-p} = f
    $$
    In particular, $ p_{k+1} / p_{k}=p/(1-p)=f$ for any
  two adjoint realizations $y_{k}, y_{k+1}$, and $k=0,2,\ldots,2^n-2$. The probabilities $p(\cdot)$ of the concatenated vectors $(0,1,s_{ij}), (1,0,s_{ij})$, etc., are
    $$
     \frac{p(0,1,s_{ij})}{p(1,0,s_{ij})} = \frac{p(0,s_{ij},1)}{p(1,s_{ij},0)} =\frac{p(s_{ij},0,1)}{p(s_{ij},1,0)} = \frac{ p^{i+1}(1-p)^{j+1} }{ p^{i+1}(1-p)^{j+1}} = 1
    $$
\item
For $n=1,2,\ldots$, $p_0 = p_0(n) =(1-p)^n$ is the probability that only $H_0$ is chosen, and $p_k = p_0 \cdot f^{\#1}$ for $k=0,\ldots,2^n-1$, where, again, $\#1$ is the number of ones in the binary representation of $k$. This means, that the vector of probabilities ${\bf p}_n = (p_0,p_1,\ldots,p_{2^n -1})$ can be written as follows:
    \begin{eqnarray*}
    {\bf p}_n &=& p_0 \cdot (1;f;f,f^2;f,f^2,f^2,f^3;f,f^2,f^2,f^3,f^2,f^3,f^3,f^4;\ldots; \\
    && f,f^2,f^2,f^3,\ldots,f^{n-2},f^{n-1},f^{n-1},f^n) = p_0 \cdot {\bf f}_n,
     \end{eqnarray*}

\item More explicitly, with ${\bf p}_0=1$, the vector ${\bf f}_n$ has dimension $2^n$ and obeys the recursive relation ${\bf f}_0=1$, and ${\bf f}_n = ({\bf f}_{n-1},f \cdot {\bf f}_{n-1})$ for $n=1,2,\ldots$ Thus its components can be calculated with the help of the following scheme, which may be interpreted as a geometric version of Pascal's triangle.\footnote{Pascal named his triangle ``triangle arithmetique'', thus, at least in French, it is straightforward to name the above multiplicative structure ``triangle geometrique''. Since row $n$ has $2^n$ entries, the geometric triangle is a true triangle on the ld scale.}
    $$
    \begin{array}{lccccccccccccccc}
   n=0: &&&&& &&& 1  &  & &&&&& \\
   n=1: &&&& 1 &&&& | &&&& f &&&\\
   n=2: && 1 && | && f && || && f && | && f^2 & \\
   n=3: & 1 & | & f & || & f & | & f^2 & ||| & f & | & f^2 & || & f^2 & | & f^3 \\
   & \multicolumn{15}{c}{\ldots} \\
    \end{array}
    $$
    Every row has $2^n$ entries. Note that the left and the right of every $|$ are ``separated'' by the factor $f$ in the following sense: First $\left[ | \right]$, $1/f=f/f^2= f^2 / f^3 =\ldots$, or, equivalently, $1 \cdot f = f; f \cdot f = f^2; f^2 \cdot f = f^3$, etc. Second $\left[||\right]$, $(1,f) \cdot f = (f, f^2); (f,f^2) \cdot f = (f^2,f^3),(f^2,f^3)\cdot f = (f^3,f^4)=\ldots$, etc. Third $\left[|||\right]$, $(1,f,f,f^2) \cdot f = (f,f^2,f^2,f^3);  (f,f^2,f^2,f^3) \cdot f = (f^2,f^3,f^3,f^4);$ etc.

    \item One may construct successive rows of (iii) in a rather elementary way: Start with a single 1 in the very first row. Then, fork every entry of row $n$ into two, by multiplying each entry with $1$ and $f$, upon moving to the left or to the right, respectively. It is quite remarkable that this ``local'' (stochastic) view is equivalent to the ``global'' (weaving) view taken in the definition.\footnote{It may be noted that the ``weaver'' is similar to the ``baker'' in dynamic system theory. In particular, in both cases a locally defined transformation is closely related to global patterns. Theorem \ref{mandelbrot} connects the stochastic and the dynamic points of view.}

\item Applying the logarithm to the base $f$ to every entry of the geometric triangle yields the exponents, i.e., the following numbers:
    $$
    \begin{array}{l|cccccccccccccccccr}
  n && & &&&&  & &  &  &  & &&&&& & \mbox{Sum}\;\; s_n  \\\hline
   0 &&& &&&&  & &  & 0 &  & &&&&& & 0 \\
  1 & && &  & & 0 &&&& | &&&& 1 && & & 1\\
 2 &&&& 0 && | && 1 && || && 1 && | && 2 && 4\\
 3 & &&  0 & | & 1 & || & 1 & | & 2 & ||| & 1 & | & 2 & || & 2 & | & 3 & 12 \\
    \multicolumn{19}{c}{\ldots} \\
    \end{array}
    $$

In general, $s_0=0$, and $s_{n+1}=2 s_n +2^n$ for $n=0,1,\ldots$ That is, one obtains the sequence $0,1,4,12,32,80,192,448,1024,2304,\ldots$
\end{enumerate}
\end{thm}

Proof: (i) is proven in the statement of the theorem. However, (i) is also obvious, since the positions of the numbers 0 and 1 are irrelevant for the probabilities in question. In particular, for $k=0,2,\ldots,2^n-2$, the binary representations of $k$ and $k+1$ differ in exactly one position.

(ii) Using theorem \ref{weaver} (ii), one obtains immediately

$p_k=p^{\#1} (1-p)^{\#0} = p^{\#1} (1-p)^{n-(\#1)} = (1-p)^n \frac{p^{\#1}}{(1-p)^{\#1}} = p_0 f^{\#1}$

(iii) is a consequence of self-similarity. Since the binary representations of 0 and $2^{n-1}$, but also of 1 and $2^{n-1}+1$, etc., only differ by a single one,
\begin{eqnarray*}
  {\bf p}_{n} &=& (p_0,\ldots,p_{2^{n-1}-1};p_{2^{n-1}},\ldots,p_{2^{n}-1}) = (p_0,\ldots,p_{2^{n-1}-1};f p_0, f p_1 , \ldots,f p_{2^{n-1}-1}) \\
  &=&  ({\bf p}_{n-1}, f {\bf p}_{n-1})  = (p_0{\bf f}_{n-1}, f p_0 {\bf f}_{n-1} ) = p_0 ({\bf f}_{n-1}, f {\bf f}_{n-1} )
\end{eqnarray*}
Since, again by (ii), also ${\bf p}_{n} = p_0 {\bf f}_{n}$, the desired result follows.

One may also prove (iii) by induction on $n$: First, $p_1 = f p_0$, and thus $(p_0,p_1)=(p_0,f p_0) = p_0 (1,f)$.
 Second, the binary representation of any $k \in \{0,\ldots,2^n-1\}$ is a vector ${\bf b}_n = (b_{n-1},\ldots,b_0)$. Let $\#1$ be the number of ones in ${\bf b}_n$. With probability $1-p$, the next selection leads to $(0,{\bf b}_n)$, and with probability $p$ this selection results in $(1,{\bf b}_n)$. Since in the first case, the number of ones does not change, and in the second case, the number of ones increases by one, we obtain on the one hand (to the left), $p_{i,n+1}  = p_{0,n+1} f^{\#1} = (1-p)^{n+1} f^{\#1} = (1-p) p_{0,n} f^{\#1} =(1-p) p_{i,n}$ for $0 \le i \le 2^n-1$. This is tantamount to ${\bf f}_n$ being reproduced as the first half of ${\bf f}_{n+1}$. (Upon moving from $n$ to $n+1$ the exponent of $f$ does not change.) On the other hand (to the right), $p_{i,n+1}  = p_{0,n+1} f^{(\#1) +1} =  (1-p)^{n+1} f^{\#1} p/(1-p)= p(1-p)^n f^{\#1} =p p_{0,n} f^{\#1} = p p_{i,n}$ for $2^n \le i \le 2^{n+1}-1$. The additional factor $f$ means that the second half of ${\bf f}_{n+1}$ has to be $f \cdot{\bf f}_n$.

(iv) The proof is by induction on $n$. For $n=0$ there is nothing to prove, and the equivalence is obvious for $n=1$. By the inductive assumption, the vector occurring on line $n$, having length $2^n$, has the form ${\bf w}_n =({\bf l}_{n-1},{\bf r}_{n-1})=({\bf l}_{n-1}, f \cdot {\bf l}_{n-1})$ where ${\bf l}_{n-1}$ is a vector of length $2^{n-1}$. In other words, $r_k/l_k = f$ for $k=1,\ldots,2^{n-1}$.

    Local splits (see the definition given in the statement of the theorem) produce a vector ${\bf w}_{n+1}$ of length $2^{n+1}$. 
    Since, locally, a step to the left reproduces the numbers, and a step to the right multiplies any two entries on tier $n$ with the same factor $f$, we also have, because of the inductive assumption, $w_{2^n+k}/w_{k}=f$ for $k=1,\ldots,2^{n}$. 
    Therefore ${\bf w}_{n+1} = ({\bf l}_n,f \cdot {\bf l}_n)$.

(v) Straightforward induction on $n$ yields the recursive formula. $\diamondsuit$

\vspace{3ex}
\begin{thm}\label{weaver2}

(Further properties of the weaver's distribution).

Given the assumptions and the notation of theorem \ref{weaver}, one obtains
\begin{enumerate}
    \item The probabilities corresponding to row $n$ can be constructed by the following simple scheme:
     $$
     \hspace{-15ex}
    \begin{array}{ccccccccccccccc}
    &&&&  & &  & 1 &  & &&&&& \\
    &&& 1-p &&&& | &&&& p & &&\\
    & \text{\begin{small}$(1-p)^2$\end{small}} && | && \text{\begin{small}$p(1-p)$\end{small}}  && || && \text{\begin{small}$p(1-p)$\end{small}}  && | && \text{\begin{small}$p^2$\end{small}}  & \\
\text{\begin{tiny}$(1-p)^3$\end{tiny}}   & | & \text{\begin{tiny}$p(1-p)^2$\end{tiny}}   & || & \text{\begin{tiny}$p(1-p)^2$\end{tiny}}
 & | & \text{\begin{tiny}$p^2(1-p)$\end{tiny}}  & ||| & \text{\begin{tiny}$p(1-p)^2$\end{tiny}} & | & \text{\begin{tiny}$p^2(1-p)$\end{tiny}}
 & || & \text{\begin{tiny}$p^2(1-p)$\end{tiny}}  & | & \text{\begin{tiny}$p^3$\end{tiny}}  \\
    \multicolumn{15}{c}{\ldots} \\
    \end{array}
    $$
 Global interpretation: ${\bf p}_{n+1} = ((1-p){\bf p}_n,p {\bf p}_n)$. Local interpretation: Start with mass 1 in the very first (the zeroth) row. Then, fork every probability of row $n$ into two, by multiplying each entry with $1-p$ (to the left) and $p$ (to the right), respectively.

    \item For $p> 1/2,$ the sequence $p_0,fp_0,f^2p_0,\ldots$ increases. Accordingly, for $p< 1/2,$ we have $f<1$. Therefore the sequence $p_0,fp_0,f^2p_0,\ldots$ decreases. If $p=1/2$, all probabilities coincide, i.e. we obtain the discrete uniform distribution on the values $y_k=k/(2^n-1)$; $p_k=1/2^n$ for $k=0,1,\ldots,2^n-1$.
 \item If $p> 1/2,$ the modus occurs in one, and the median is larger than $1/2$. Vice versa, if $p< 1/2,$ the modus occurs in zero, and the median is less than $1/2$.
\item Symmetry: Suppose $Y \sim W(n,p)$ and $Y^{'} \sim W(n,1-p)$. Then $P(Y=y_k) = P(Y^{'}=y_{2^n-1-k})$ for $k=0,\ldots,2^n -1$. 

    \item Distribution function of $W(n,p)$: For all $n\ge2$ and $k=0,\ldots,2^n$ define $v_{k,n}=k/2^n$. For every fixed $n$, the mass left and right to $v_{k,n}$ is constant for every $m \ge n$, and so is the value of $F(v_{k,n})$. In particular, $F(v_{1,1})=F(1/2)=(1-p)$ for all $n \ge 1$; $F(v_{1,2})=F(1/4)=(1-p)^2$, $F(v_{3,2})=F(3/4)=1-p^2$ for all $n\ge 2$; $F(v_{1,3})=F(1/8)=(1-p)^3$; $F(v_{3,3})=F(3/8)=(1-p)^2+p(1-p)^2$, $F(v_{5,3})=F(5/8)=(1-p)+p(1-p)^2$, $F(v_{7,3})=F(7/8)=1-p^3$ for all $n\ge 3$, etc.
        \item For $k=0,\ldots,2^n-1$, the total mass in every interval $[v_{k,n},v_{k+1,n}]$ remains constant. Moreover, the mass in this interval is located at a single point, $y_k=k/(2^n-1)$.
            \item Distribution of the jumps (stick heights): $F_n$ has $2^n$ points of discontinuity. If $p=1/2$ there is a constant jump height $h=1/2^n$. Otherwise, there are $n+1$ different jumps sizes, actually $h_j=p^j(1-p)^{n-j}$ for $j=0,\ldots,n$, having a binomial distribution. That is, there is 1 jump of size $h_0=(1-p)^n$, there are $\binom{n}{1}=n$ jumps of size $h_1=(1-p)^{n-1}p$, etc.
      \end{enumerate}
\end{thm}

Proof: (i) For $n=1,2,\ldots$, we have $p_0=p_0(n)=(1-p)^n$ for the leftmost probability (only $H_0$ is selected). Applying the geometric triangle yields the result. (ii) we have $p<1/2 \Rightarrow f>1$. Thus the mass in $y_1$ exceeds the mass in $y_0=0$ by the factor $f$, and the result follows straightforwardly.
(iii) is due to self-similarity. The claim for the modus can also be shown directly, since, if $p< 1/2,$ we have $(1-p)^n < (1-p)^{n-k} p^k < p^n$.
(iv). Exchanging the roles of zeros and ones, and replacing $p$ by $1-p$ yields the same distribution. In other words: The mirror image of $W(n,p)$ by the symmetry axis $y=1/2$ is $W(n,1-p)$. (v) follows immediately from the geometric triangle. Geometrically speaking, the unit interval on the abscissa is successively halved. At the same time, the unit interval on the ordinate is successively split according to the ratio $f$. Thus, for finite $n \ge 1$, one obtains a step function with $2^n$ jumps.
(vi) holds because of the local interpretation of the geometric triangle: Each split can be interpreted as distributing the mass $p_k$ in $y_k$ to the points $y_{2k,n+1}$ and $y_{2k+1,n+1}$ in that same interval. Graphically, the stick of height $p_k$ in $y_{k,n}$ is broken into two sticks of heights $(1-p)p_k$ and $p \cdot p_k$, located in $y_{2k,n+1}$ and $y_{2k+1,n+1}$, respectively. (vii) is due to construction.
$\diamondsuit$

 In the last theorem, (i) is  the classical ``binomial'' or ``p modell'' first described by \citet{de51, de53}. Note, however, that although there is a close relationship with the binomial distribution, a weaver's distribution is based on progressive sampling, and thus obtains $2^{n}$ values. A closer look reveals that there are two scales involved, the first given by discrete ``time'', i.e. the number of observations $t=2^n$, and the second by logarithmic time, that is, $\ld 2^n =n$, the number of distribution-selects.

\section{Expected value}

\begin{thm}\label{expectedvalue}

Let $Y_n \sim W(n,p)$. Then, for every $n \ge 1$, the expected value of $Y_n$ is $p$.
\end{thm}

Proof: Let $\mu =EY_n$. One may decompose $\mu$ into a sum of $n$ terms $t_0,\ldots,t_{n-1}$, where the index $j$ counts the number of zeros in the corresponding binary vector ${\bf b}_n=(b_{n-1},\ldots,b_0)$, that is, $j=n-\sum_{i=0}^{n-1} b_i$. More precisely, $\mu = \sum_{j=0}^{n-1} t_j =\sum_{j=0}^{n-1} p_j \cdot y_{[j]}$ where
$y_{[j]}$ is the sum of all realizations with corresponding probability mass $p_j$.

$j=0$: There is only one vector of dimension $n$ without the entry zero, i.e., ${\bf b}_n=(1,\ldots,1)$. The corresponding probability is $p^n$ and thus $t_0=1 \cdot p^n$

$j=1$: We have to consider the sum of all realizations of $Y$ that occur with probability $p_1=p^{n-1}(1-p)$, i.e. all binary sequences of length $n$, having exactly one zero. Thus
\begin{eqnarray*}
y_{[1]} &=& \left( 2^n -1 -2^0 + 2^n -1 -2^1 + 2^n-1 -2^2 + \ldots + 2^n-1-2^{n-1} \right) / (2^n-1) \\
&=& (n2^n-n-\sum_{i=0}^{n-1} 2^i )/ (2^n-1) =  (n (2^n -1)-(2^n-1)) / (2^n-1) = n-1
\end{eqnarray*}
More intuitively, the number $2^n-1$ is represented by a vector of $n$ successive ones in the binary system. In the last equation we are looking for all sequences of length $n$ with exactly one zero. There are exactly $n$ such sequences, with the zero placed in every possible position. Thus their sum is $n(2^n-1)-(2^n-1)=(n-1)(2^n-1)$. Dividing by $2^n-1$ yields the result, and $t_1 = (n-1) p^{n-1}(1-p)$.

$j=2$: There are $\binom{n}{2}$ ways to place exactly two zeros in a string of length $n$. Without the zeros, the sum of these sequences would be $\binom{n}{2} (2^n-1)$. However, for every ``chain'' of zeros we have to subtract $\sum_{i=0}^{n-1} 2^i$, and there are $\frac{2}{n}\cdot\binom{n}{2}$ such chains. Thus
\begin{eqnarray*}
y_{[2]} &=& \left(\binom{n}{2} (2^n -1) - \frac{2}{n} \binom{n}{2} \sum_{i=0}^{n-1} 2^i \right) / (2^n-1) \\
&=& \left(\frac{n(n-1)}{2} (2^n -1) - (n-1) (2^n-1) \right) / (2^n-1) \\
&=& (n-1) \left(\frac{n}{2}-1 \right) =  \binom{n-1}{2} \\
\end{eqnarray*}

Thus $t_2=p_2 y_{[2]} = \binom{n-1}{2} p^{n-2}(1-p)^2.$

$j=3$: There are $\binom{n}{3}$ ways to place exactly three zeros in a string of length $n$. Without the zeros, the sum of these sequences would be $\binom{n}{3} (2^n-1)$. However, for every ``chain'' of zeros we have to subtract $\sum_{i=0}^{n-1} 2^i=2^n-1$, and there are $\frac{3}{n}\cdot\binom{n}{3}$ such chains. Thus
\begin{eqnarray*}
y_{[3]} &=& \left(\binom{n}{3} (2^n -1) - \frac{3}{n} \binom{n}{3} (2^n-1) \right) / (2^n-1)
= \binom{n}{3}  - \binom{n-1}{2}  =  \binom{n-1}{3}, \\
\end{eqnarray*}
and therefore $t_3=\binom{n-1}{3} p^{n-3}(1-p)^3.$

With exactly the same argument, we find all further terms $t_j$, and finally
$$
t_{n-3}= \binom{n-1}{n-3} p^{3}(1-p)^{n-3}, t_{n-2}= \binom{n-1}{n-2} p^{2}(1-p)^{n-2}, t_{n-1}=  p(1-p)^{n-1}
$$
Putting everything together, we get with the help of the Binomial theorem:
\begin{eqnarray*}
\mu &=& \sum_{j=0}^{n-1} t_j = p^n + \sum_{j=1}^{n-1} \binom{n-1}{j} p^{n-j} (1-p)^j = p^n + \sum_{j=0}^{n-1} \binom{n-1}{j} p^{n-j} (1-p)^j -p^n \\
&=& p \sum_{j=0}^{n-1} \binom{n-1}{j} p^{(n-1)-j} (1-p)^j =p \;\;\;\;\;\; \diamondsuit
\end{eqnarray*}

\section{Variance}

After the first step, the distribution of the conditional expected values is $B(p)$. For any random variable $X$ with values in the unit interval, and $EX=p$, this distribution has maximum variance $p(1-p)$. Upon weaving, probability mass is successively shifted into the unit interval, and thus variance decreases. On the other hand, every bifurcation may increase the variance term. However, both effects combined could result in an (net) monotone decrease of variance up to a certain point. (For concrete values, see the table p. \pageref{concrete}.) Moreover, there should be a limit variance $\sigma^2= c p (1-p)$ with $c <1$.

\begin{thm}\label{variancetheorem}
Let $Y_n \sim W(n,p)$. Then the variance of this random variable is
\begin{equation}\label{variance}\sigma^2(Y_n) =\frac{\sum_{i=0}^{n-1} 2^{2i}}{(2^n-1)^2} p (1-p)\end{equation}
\end{thm}

Proof: If we interpret $k=\sum_{i=0}^{n-1} b_i$ as a binary number of length $n$, the $i+1$th step of the above selection scheme defines its $i$th digit (from the right to the left, $i=0,\ldots,n-1$). Since, by construction, the digits are independent, every step contributes a certain amount to the overall variance, independent of all the other steps. This means that the total variance can be decomposed into $n$ parts $\sigma_0^2,\ldots,\sigma_{n-1}^2$ that accrue to the total variance. The variance contributed by the $i$th digit is the difference between $(?\cdots?1?\cdots?)$ and $(?\cdots?0?\cdots?)$, where the question marks denote arbitrary other binary digits (the same for both numbers). 

As a typical example, consider the case $n=3$. The first step introduces variance that can be assessed by means of considering two adjoint realizations, e.g., the values $0=(000)_2$ and $1/7=(001)_2 / (111)_2$. This results in
$$
\sigma_0^2 = p \left( \frac{1}{7}-\frac{1}{7} p \right)^2 + (1-p) \left( 0-\frac{1}{7} p \right)^2 = \frac{1}{49} p(1-p)  = \left( \frac{1}{7} \right)^2 p (1-p)
$$
By the same token, the variance produced by the second step can be measured by two realizations that only differ in the second component of their binary representation, e.g., the values $0=(000)_2$ and $2/7=(010)_2 / (111)_2$. This gives
$$
\sigma_1^2 = p \left( \frac{2}{7}-\frac{2}{7} p \right)^2 + (1-p) \left( 0-\frac{2}{7} p \right)^2 = \frac{4}{49} p(1-p)  = \left( \frac{2}{7} \right)^2 p (1-p)
$$
Finally, since the variance produced by the last step (4 bifurcations) is the same for all their descendants, it suffices to consider just one of these forks, e.g., the values $0=(000)_2$ and $4/7=(100)_2 / (111)_2$. This leads to
$$
\sigma_2^2 = p \left( \frac{4}{7}-\frac{4}{7} p \right)^2 + (1-p) \left( 0-\frac{4}{7} p \right)^2 = \frac{16}{49} p(1-p) = \left( \frac{4}{7} \right)^2 p (1-p)
$$
Putting everything together, we obtain
$$
\sigma^2(Y_3) = \sigma_0^2 + \sigma_1^2 + \sigma_2^2 = (1+4+16) p(1-p) / 7^2
$$
Therefore, in general, $\sigma^2(Y_n) = \sum_{i=0}^{n-1} \sigma_i^2$, where $$\sigma_i^2 = p \left( \frac{2^i}{2^n-1}-\frac{2^i}{2^n-1} p \right)^2 + (1-p) \left( 0-\frac{2^i}{2^n-1} p \right)^2  =(2^i)^2  p(1-p) / (2^n -1)^2.$$
$\diamondsuit$

\vspace{2ex}
Note that the numerator shows an additive analog to factorials: For factorials, $n! = (n-1)! \cdot n$ holds. For the numerator, we have $\sum_{i=0}^{n} (2^{i})^2 = \sum_{i=0}^{n-1} 2^{2i} + 2^{2n}$.

\begin{cor}\label{moments}
$EY_n^2$ exists, and so do all higher moments $EY_n^j$ for $j \ge 1$.
\end{cor}

Proof: For fixed $n$, all realizations $y_k$ are in the unit interval. Thus $y_k \ge y_k^2 \ge y_k^3 \ge \dots$, with strict inequality if $0 < y_k < 1$. Therefore $0 < EY_n^i < EY_n^j$ if $i > j$. $\diamondsuit$

\begin{lem}\label{limitvariance}
The limit of the variance term is $\frac{1}{3}p(1-p)$
\end{lem}

Proof: Considered as a function of $n$, $\sigma^2(Y_n)$ is monotonically decreasing. Since it is also nonnegative, it is clearly convergent. Moreover, a straightforward induction on $n$ shows that $\sum_{i=0}^{n-1}2^{2i} = (2^{2n}-1)/3$, thus
\begin{eqnarray*}
  \frac{\sigma^2(Y_n)}{p (1-p)} &=& \frac{\sum_{i=0}^{n-1}2^{2i} }{ (2^n -1)^2}= \frac{2^{2n}-1}{3 (2^{2n}-2^{n+1}+1)}
   = \frac{(2^{2n}-1) / 2^{2n}}{3 (2^{2n}-2^{n+1}+1) / 2^{2n}} ,  \\
\end{eqnarray*}
which converges to $1/3$ if $n \rightarrow \infty$. $\diamondsuit$

\section{Limit distribution}
   Since, owing to theorem \ref{weaver2}, the distribution function $F_n$ is well-known for all values $v(k,n)$, it is easy to pass to the limit. The limit function $F$ obviously is a distribution function.

\begin{thm}\label{limitdistribution} (The weaver's hem)

Let $Y$ be the limit of $(Y_n)$, defined by its distribution function $F=\lim_{n \rightarrow \infty} F_n$. For obvious reasons, the corresponding distribution, i.e., $Y \sim W(p)$, may be called ``the weaver's hem.'' Its first moments are $EY=p$ and $\sigma^2(Y)=p(1-p)/3$. Except for the case $p=1/2$, when the discrete uniform distribution becomes the continuous uniform distribution on the unit interval (and thus $F$ is the identity function there), $F$ has no density with respect to Lebesgue measure.
\end{thm}

Proof:
Because of $EX=\int_0^1 (1-G(x)) dx$ for any distribution function $G$ on the unit interval, and $F_n\rightarrow F$, we also have $EY=p$ for the weaver's hem. An analogous argument for the second moment and Theorem 7 yields $\sigma^2(Y)=p(1-p)/3.$

   Heuristically, if $p > 1/2$, consider the interval $[0,1/2[$. The mass $1-p$ available there is shifted to the left. Thus the distribution function grows rapidly first, but hardly grows near $1/2$. Now consider the interval $]1,2,1]$. There, mass $p$ is available and systematically shifted to the left. Thus the distribution function grows rapidly near $1/2$, but very slowly in the vicinity of $1$. Thus the distribution function has a salient point in $1/2$ and cannot be differentiated there. The same holds for all $v(k,n)$. Since the set of these points lies dense in the unit interval, there should be no density.

   Formally, consider an interval of length $1/(2^n-1)$ about every $y_k=k/(2^n-1)$. The probabilities at the latter points all have the form $p^j (1-p)^{n-j}$, where $0 \le j \le n$. Thus the density is given by\footnote{Note, that in the classical De Moivre-Laplace theorem, one has to consider $\binom{n}{j}p^j (1-p)^{n-j}$.}
   $$
   f_n(j)= (2^n-1) p^j (1-p)^{n-j}.
   $$
   If $p=1/2$, $f_n(j)=1-1/2^n$. Thus there is a finite limit which does not depend on $j$. Since $p^j (1-p)^{n-j} \rightarrow 0$ for all $j$ if $n \rightarrow \infty$, it suffices to consider $2^n p^j (1-p)^{n-j}$. W.l.o.g. let $p<1/2$. Then we obtain for the interval about one $2^n p^n =(2p)^n \rightarrow 0$, and for the interval about zero we get $2^n(1-p)^n=(2-2p)^n \rightarrow \infty$ for $n \rightarrow \infty$. Similarly, the densities for all the intervals in between either converge to zero or diverge. Thus, no limit density exists if $p \neq 1/2$. $\diamondsuit$

\begin{thm}\label{mandelbrot}
   The limit distribution $W(p)$ is equivalent to Mandelbrot's well-known ``binomial measure'' (a multifractal).
\end{thm}

Proof: Mandelbrot's measure is defined as the limit of the following iterative process (cf. \citet{ma74}, p. 329): For $n=0$, start with the uniform distribution on the unit interval. Next, the proportion $1-p$ is uniformly distributed on the interval $(0,1/2)$, and the proportion $p$ is uniformly distributed on the interval $(1/2,1)$. Then, one splits the masses further (locally) according to the geometric triangle, i.e., mass $(1-p)^2$ to the interval $(0,1/4)$, mass $(1-p)p$ to the interval $(1/4,1/2)$, mass $p(1-p)$ to the interval $(1/2,3/4)$, and mass $p^2$ to the interval $(3/4,1)$, etc.

Thus, Mandelbrot's scheme and weaving, if interpreted locally (see theorem \ref{weaver2}, (v) and (vi)), refer the same mass to every interval $[v_{k,n};v_{k+1,n}]$. Since these intervals shrink to zero, the limit distributions have to coincide. It is well-known that the corresponding distribution function, with the exception of $p=1/2$, has no density (cf. \citet{ri99}).
$\diamondsuit$

\section{The complete process}

So far, we have mainly considered the distribution of the (conditional) expected values, $Y_n = E({\bar X}_n | B)$, or, equivalently, the case of two one-point distributions located in $\mu_0$ and $\mu_1$. Looking at ${\bar X}_n$, however, there is not just variance between the populations $H_0$ and $H_1$, that we have considered so far, but also within each of these populations, $\sigma^2(H_0) = \sigma_0^2$ and $\sigma^2(H_1)=\sigma_1^2$, say, contributing to the total variance.

In complete generality, i.e., without specific distributional assumptions or any particular sampling scheme, suppose $n_0$ independent observations $Z_1,\ldots,Z_{n_0}$ come from the first population, and $n_1$ independent observations $Z_1^{'},\ldots,Z_{n_1}^{'}$ come from the second population, $n=n_0+n_1$.
At this point of sampling, the combined distribution is a mixture $M$ given weight $n_0/n$ to the sample from $H_0$, and weight $n_1/n$ to the sample from $H_1$. In particular,
$$
{\bar X}_n = \frac{\sum_{i=1}^n X_i}{n} = \frac{\sum_{i=1}^{n_0} Z_i +\sum_{i=1}^{n_1} Z^{'}_i}{n}
= \frac{n_0}{n}\frac{\sum_{i=1}^{n_0} Z_i}{n_0} +\frac{n_1}{n}\frac{\sum_{i=1}^{n_1} Z^{'}_i}{n_1}
$$
Thus we get the expected value (total mean)
\begin{eqnarray}\label{mixedmean}
\mu &=& E{\bar X}_n = E[E({\bar X}_n | M)] = \frac{n_0}{n}\mu_0 + \frac{n_1}{n}\mu_1
\end{eqnarray}

and variance
\begin{eqnarray}\label{mixedvariance}
\sigma_n^2 &=& \sigma^2 (E({\bar X}_n |M)) + E(\sigma^2({\bar X}_n | M)) \nonumber \\
 &=& \frac{n_0}{n} (\mu_0 -\mu)^2 + \frac{n_1}{n} (\mu_1 -\mu)^2 + \frac{n_0}{n} \frac{\sigma_0^2}{n_0} + \frac{n_1}{n} \frac{\sigma_1^2}{n_1}  \\\nonumber
\end{eqnarray}

\begin{thm}\label{amproperties} (Expected value and variance)

  With the assumptions of theorem \ref{weaver} , $E{\bar X_n} = \mu =p$
  and
  \begin{equation}\label{varianz}
  \sigma^2({\bar X}_n) = p (1-p)+ \frac{\sigma_0^2 + \sigma_1^2}{2^n-1}
 \end{equation}
\end{thm}

Proof:
Given progressive sampling, after $n$ selections, there are $2^n$ mixed distributions $Q_k$ with the proportion $\lambda_k =k/(2^n-1)=y_k$ of observations coming from $H_1$  ($k=0,\ldots,2^n-1$). Distribution $Q_k$ occurs with probability $p_k$, where $p_k$ comes from a $W(n,p)$ distribution. If $Z_k \sim Q_k$, and $\mu_k=EZ_k$, equation $(\ref{mixedmean})$ translates into
\begin{equation}\label{erwartung}
\mu = \sum_{k=0}^{2^n-1} p_k \mu_k = \sum_{j=0}^{n-1} p_j y_{[j]} =p
\end{equation}
due to theorem \ref{expectedvalue}, where we used the notation of the latter theorem, that is, $y_{[j]}$ is the sum of all $\mu_k = y_k= E({\bar X}_n | B_n=(b_{n-1},\ldots,b_0))$ with corresponding probability mass $p_j$. In other words, the sum extends over all vectors $(b_{n-1},\ldots,b_0)$ containing exactly $j$ zeros, $j=n-\sum_{i=0}^{n-1} b_i$.

The first part of equation $(\ref{mixedvariance})$, capturing the variance between the $Z_k$, reads
$$
\sigma^2 (E({\bar X}_n |B_n)) = \sum_{k=0}^{2^n-1} p_k (\mu_k - \mu)^2 =\frac{\sum_{i=0}^{n-1} 2^{2i}}{(2^n-1)^2} p (1-p)
$$
due to theorem \ref{variancetheorem}.
Finally, the second part of equation $(\ref{mixedvariance})$, accounting for the variance within the mixtures, becomes
\begin{eqnarray*}
E(\sigma^2({\bar X}_n | B_n)) &=& \sum_{k=0}^{2^n-1} p_k \sigma^2({\bar X}_n | B_n=(b_{n-1},\ldots,b_0))) \\
\end{eqnarray*}
For every fixed $k=(b_{n-1},\ldots,b_0)_2$, $Q_k$ is a mixture with $k=\sum_{i=0}^{n-1} b_i 2^{i}$ observations from $H_1$. Using $\mu=0$ and $\mu=1$, equation $(\ref{mixedmean})$ simplifies to $\mu_k=\lambda_k = k/(2^n-1)$ and the variance of $Z_k$, again according to equation $(\ref{mixedvariance})$, is
\begin{eqnarray*}
\sigma^2({\bar X}_n | B_n=(b_{n-1},\ldots,b_0)) & =&   (1-\lambda_k) \lambda_k^2 + \lambda_k (1-\lambda_k)^2  +
 (1-\lambda_k) \frac{\sigma_0^2}{2^n-1-k} + \lambda_k \frac{\sigma_1^2}{k} \\
 & =&   (1-\lambda_k) \lambda_k + \frac{\sigma_0^2}{2^n-1} +  \frac{\sigma_1^2}{2^n-1}
\end{eqnarray*}

Altogether we obtain the preliminary result
\begin{equation}\label{preliminary}
\sigma^2 ({\bar X}_n) =\frac{\sum_{i=0}^{n-1} 2^{2i}}{(2^n-1)^2} p (1-p) + \sum_{k=0}^{2^n-1} p_k  \left( \lambda_k(1-\lambda_k)+ \frac{\sigma_0^2 + \sigma_1^2}{2^n-1}  \right)
\end{equation}

Now
\begin{eqnarray*}
\sum_{k=0}^{2^n-1} p_k \lambda_k (1-\lambda_k) &=& \sum_{k=0}^{2^n-1} p^{\sum_{i=0}^{n-1} b_i} (1-p)^{n-\sum_{i=0}^{n-1} b_i}
\frac{k}{2^n-1}\frac{2^n-1-k}{2^n-1} \\
&=& \frac{1}{(2^n-1)^2}\sum_{k=1}^{2^n-2} p^{\sum_{i=0}^{n-1} b_i} (1-p)^{n-\sum_{i=0}^{n-1} b_i}
\left( \sum_{i=0}^{n-1} b_i 2^{i} \right) \left(  \sum_{i=0}^{n-1} (2^n -1 - b_i 2^{i}) \right) \\
&=& \frac{p(1-p)}{(2^n-1)^2} \left\{  1 \cdot (2^n-2) (1-p)^{n-2} + 2 \cdot (2^n-3) (1-p)^{n-2}  \right. \\
&& + \left. 4 \cdot (2^n-3) p (1-p)^{n-3} + \ldots  +(2^n-3) \cdot 2 \cdot p^{n-2} + (2^n-2) \cdot 1 \cdot p^{n-2}  \right\} \\
\end{eqnarray*}
The last term in brackets can be rearranged:
\begin{eqnarray*}
\left\{ \ldots \right\} &=&  (1-p)^{n-2} [1 \cdot (2^n-2) + 2 \cdot (2^n-3)+ 4 \cdot (2^n-5) + \ldots + 2^{n-1} (2^{n-1}-1)]  \\
&& + p (1-p)^{n-3} [3 \cdot (2^n-4) + 5 \cdot (2^n-6) +\ldots + (2^{n-1}+2^{n-2}) \cdot (2^n-1 -(2^{n-1}+2^{n-2})) ]   \\
&& + \ldots +  p^{n-2} [(2^{n-1}-1) 2^{n-1} + \ldots + (2^{n}-5) \cdot 4 + (2^n-3) \cdot 2+ (2^n-2) \cdot 1] \\
&=& \binom{n-2}{0} (1-p)^{n-2} \left( \sum_{j=0}^{n-1} 2^j (2^n-1-2^j) \right) + \binom{n-2}{1} p (1-p)^{n-3} \left( \sum_{j=0}^{n-1} 2^j (2^n-1-2^j) \right) \\
&& + \ldots \binom{n-2}{n-3} p^{n-3} (1-p) \left( \sum_{j=0}^{n-1} 2^j (2^n-1-2^j) \right) + \binom{n-2}{n-2} p^{n-2}\left( \sum_{j=0}^{n-1} 2^j (2^n-1-2^j) \right) \\
&=& \sum_{j=0}^{n-1} 2^j (2^n-1-2^j) (1-p + p)^{n-2} =\sum_{j=0}^{n-1} 2^j (2^n-1-2^j)
\end{eqnarray*}
so that
$$
\sum_{k=0}^{2^n-1} p_k \lambda_k (1-\lambda_k) = p(1-p) \sum_{j=0}^{n-1} 2^j (2^n-1-2^j) / (2^n-1)^2
$$
and equation (\ref{preliminary}) becomes
\begin{eqnarray}\label{weavemix}\nonumber
\sigma_n^2 ({\bar X}_n) &=& \frac{\sum_{i=0}^{n-1} 2^{2i}}{(2^n-1)^2} p (1-p) + \sum_{k=0}^{2^n-1} p_k   \lambda_k(1-\lambda_k)+  \sum_{k=0}^{2^n-1} p_k \frac{\sigma_0^2 + \sigma_1^2}{2^n-1} \\
&=& \frac{\sum_{i=0}^{n-1} 2^{2i}}{(2^n-1)^2} p (1-p) + \frac{\sum_{j=0}^{n-1} 2^j (2^n-1-2^j) }{(2^n-1)^2} p (1-p)+  \frac{\sigma_0^2 + \sigma_1^2}{2^n-1} \\\nonumber
&=& p(1-p) +  \frac{\sigma_0^2 + \sigma_1^2}{2^n-1} , \\\nonumber
\end{eqnarray}
where the last equation is due to the next technical lemma.
$\diamondsuit$

\begin{lem}\label{technical}
$$
\frac{\sum_{i=0}^{n-1} 2^{2i}}{(2^n-1)^2} + \frac{\sum_{j=0}^{n-1} 2^j (2^n-1-2^j)}{(2^n-1)^2} =1
$$
\end{lem}

Proof: All one has to do is rearrange the terms:
\begin{eqnarray*}
 \sum_{i=0}^{n-1} 2^{2i} &+& \sum_{j=0}^{n-1} 2^j (2^n-1-2^j) =
2^0 + 2^2 + 2^4 + \ldots + 2^{2(n-1)}  \\
&+&  2^0 (2^n-1-2^0)+ 2^1 (2^n-1-2^1)+ 2^2 (2^n-1-2^2) +\ldots + 2^{n-1} (2^n-1-2^{n-1}) \\
&=& 2^0 + 2^0 (2^n-1-2^0) + 2^2  + 2^1 (2^n-1-2^1) +  2^4  + 2^2 (2^n-1-2^2) \\
&& +\ldots +2^{2(n-1)}  + 2^{n-1} (2^n-1-2^{n-1})  \\
  &=& (2^n-1) + 2 (2^n-1)+\ldots +2^{n-1}  (2^n-1) \\
&=& (2^n-1) (1+2+\ldots +2^{n-1} ) =(2^n-1) (2^n-1) \;\;\; \;\;\; \diamondsuit \\
\end{eqnarray*}
It may be helpful to display some concrete values
$$\label{concrete}
\begin{array}{|l|r|r|r|r|c|}
  n & 2^n-1 & (2^n-1)^2 & \sum_{i=0}^{n-1} (2^i)^2 & \sum_{i=0}^{n-1} 2^i (2^n-1-2^i) & ratios \\\hline
 1 & 1 & 1 & 1 & 0 & 1 \leftrightarrow 0 \\
 2 & 3 & 9 & 4 & 5 & 4/9=0.{\bar 4} \leftrightarrow 5/9=0.{\bar 5} \\
 3 & 7 & 49 & 21 & 28 & 0.43 \leftrightarrow 0.57 \\
 4  & 15 & 225 & 85 & 140 & 0.3{\bar 7} \leftrightarrow 0.6{\bar 2}\\
 5 & 31 & 931 & 341 & 620 & 0.35 \leftrightarrow 0.65 \\
 6 & 63 & 3969 & 1365 & 2604 & 0.34 \leftrightarrow 0.66 \\
  \hline
\end{array}
$$

With respect to weaving and merging this means that one starts $(n=1)$ with a $B(p)$ distribution having expected value $p$ and variance $p(1-p)$ on the unit interval. In the next steps, this ``available'' variance is successively distributed among weaving and mixing, since due to theorem \ref{amproperties} the latter variances add up to $p(1-p)$ for $n \ge 2$.

With $n$ increasing, theorems \ref{limitdistribution} and \ref{amproperties} govern the asymptotic behaviour. That is, the first variance component due to weaving (i.e., the first term in equation (\ref{weavemix})) decreases towards $1/3$, which has the consequence that the second component due to mixing (i.e., the second term in equation (\ref{weavemix})) has to increase to $2/3$. Moreover, since the variance within the populations (i.e., the third term in equation (\ref{weavemix})) vanishes, we obtain the following result:

\begin{thm}\label{amlimit} (Limit distribution)

Given the assumptions and the notation of theorem \ref{weaver}, if $H_0$ and $H_1$ both have finite variances, then $B(p)$ is the limit distribution of the inhomogeneous (unconditional) stochastic process $({\bar X}_n)$.
\end{thm}

Proof: If there were just one population, $H_1$, say, ${\bar X}_n$ would converge to 1 almost surely. Because $H_0$ makes ${\bar X}_n$ smaller, at least in expectation, 1 has to be the uppermost accumulation point. Since, for the same reasoning, 0 is the lowest accumulation point, $P({\bar X}_n \notin [0,1]) \rightarrow 0$.

Now, because of the last theorem, for every $n$, the process is centered in $E{\bar X}_n =p$, and its variance is given by equation $(\ref{varianz})$. Obviously, $(\sigma_0^2+\sigma_1^2) / (2^n-1)$ converges to zero. Thus we are left with a limit distribution that is restricted to the unit interval, centered in $p$ and has maximum variance $p(1-p)$. These properties imply the result.
$\diamondsuit$

Intuitively, the latter result is quite obvious: If the variance within the populations vanishes, it is just the variance between the populations that is asymptotically relevant. The latter variance is due to weaving (i.e., the selection process of the populations), and conditional mixing.

Of course, if the populations $H_0$ and $H_1$ are not too complicated, it is possible to study the process $(X_n)$ in much more detail.

\end{document}